\documentclass[a4paper,10pt]{article}

\usepackage[english]{babel}
\usepackage[utf8]{inputenc}

\usepackage{amsmath}
\usepackage{amsfonts}
\usepackage{amssymb}
\usepackage{amsthm}
\usepackage{graphicx}
\usepackage{pdfpages}

\newcommand{\mytext}[1]{ \: \textrm{#1} \: }
\newcommand{\mysetdescr}[2]{\left\{ #1 \: \left| \: #2 \right. \right\} }

\newcommand{\myN}{\mathbb{N}}
\newcommand{\myNk}[1]{\underline #1}
\newcommand\mytimes{{\times}}

\newcommand{\strG}{\sqsubseteq_\Gamma}

\newcommand\myurbild[1]{{#1}^{-1}}

\def\A{{\cal A}}
\def\G{{\cal G}}
\def\H{{\cal H}}
\def\I{{\cal I}}
\def\J{{\cal J}}
\def\L{{\cal L}}
\def\S{{\cal S}}
\def\T{{\cal T}}
\def\U{{\cal U}}

\newcommand{\mygxiv}[2]{\Gamma_{#1}(#2)}
\newcommand{\gxiv}{\mygxiv{\xi}{v}}
\newcommand{\grxiv}{\mygxiv{\rho(\xi)}{v}}
\newcommand{\grgxiv}{\mygxiv{\rho_G(\xi)}{v}}

\newcommand{\myGxi}[1]{\G(#1)}
\newcommand{\Gxi}{\myGxi{\xi}}
\newcommand{\Gzet}{\myGxi{\zeta}}

\newcommand{\Txi}{\T(\xi)}

\newcommand{\grgxixv}{\mygxiv{\rho_G(\xi)}{v}}

\newcommand{\gzev}{\mygxiv{\zeta}{v}}

\newcommand{\rhoxiv}{\rho(\xi)(v)}
\newcommand{\rhoxiw}{\rho(\xi)(w)}
\newcommand{\xiv}{\xi(v)}
\newcommand{\xiw}{\xi(w)}

\newcommand{\mf}[1]{\mathfrak{ #1 }}
\newcommand{\fa}{\mf{a}}
\newcommand{\fb}{\mf{b}}
\newcommand{\fc}{\mf{c}}
\newcommand{\fC}{\mf{C}}
\newcommand{\fD}{\mf{D}}
\newcommand{\fg}{\mf{g}}
\newcommand{\fG}{\mf{G}}
\newcommand{\fP}{\mf{P}}
\newcommand{\fS}{\mf{S}}
\newcommand{\fT}{\mf{T}}
\newcommand{\fU}{\mf{U}}

\newcommand{\fTa}{\fT_a}

\newcommand{\Nin}{N^{in}}
\newcommand{\Nout}{N^{out}}

\newcommand{\ariestr}{G \in \fD'_r$, $\xi \in \H(G,R)$, $v \in V(G)}

\newcommand{\ariegxistr}{G \in \fD'_r$, $\xi \in \H(G,R)}

\def\BP{\begin{proof}}
\def\EP{\end{proof}}

\begin{document}

\theoremstyle{plain}
\newtheorem{condition}{Condition}
\newtheorem{theorem}{Theorem}
\newtheorem{definition}{Definition}
\newtheorem{corollary}{Corollary}
\newtheorem{lemma}{Lemma}
\newtheorem{proposition}{Proposition}

\title{\bf Criteria for the less-equal-relation between partial Lov\'{a}sz-vectors of digraphs}
\author{\sc Frank a Campo}
\date{\small Seilerwall 33, D 41747 Viersen, Germany\\
{\sf acampo.frank@gmail.com}}

\maketitle

\begin{abstract}
\noindent Finite digraphs $R$ and $S$ are studied with $\# \H(G,R) \leq \# \H(G,S)$ for every finite digraph $G \in \fD'$, where $\H(G,H)$ is the set of order homomorphisms from $G$ to $H$ and $\fD'$ is a class of finite digraphs. It is shown that for several classes $\fD'$ of digraphs and $R \in \fD'$, the relation $\# \H(G,R) \leq \# \H(G,S)$ for every $G \in \fD'$ is implied by the relation $\# \S(G,R) \leq \# \S(G,S)$ for every $G \in \fD'$, where $\S(G,H)$ is the set of homomorphisms from $G$ to $H$ mapping all proper arcs of $G$ to proper arcs of $H$. Under an application-oriented regularity condition, the two relations are even equivalent. A method is developed for the rearrangement of a digraph $R$, resulting in a digraph $S$ with $\# \H(G,R) \leq \# \H(G,S)$ for every digraph $G$. The method is applied in constructing pairs of partially ordered sets $R$ and $S$ with $\# \H(P,R) \leq \# \H(P,S)$ for every partially ordered set $P$. The main part of the results holds also for undirected graphs.
\newline

\end{abstract}

\section{Introduction} \label{sec_introduction}

The systematic study of homomorphisms between directed graphs (digraphs) begun around 1940 with focus on partially ordered sets (posets). In order to unify ordinal and cardinal arithmetic, Garret Birkhoff \cite{Birkhoff_1937,Birkhoff_1942} published two articles in 1937 and 1942 in which he introduced the direct sum and product of posets $P$ and $Q$ and the homomorphism set $\H(P,Q)$, together with their (later on) usual partial order relations. By doing so, he opened the rich field of ``order arithmetic'' in which the homomorphism sets $\H(P,Q)$, often notated as $Q^P$, play the role of exponentiation. Surveys about the respective state of the art are contained in J\'{o}nsson \cite{Jonsson_1982_1} from 1982, Duffus \cite{Duffus_1984} from 1984, and McKenzie \cite{McKenzie_2003} from 2003.

The interest in homomorphisms between undirected graphs started in the sixties with the pioneering work of Sabidussi \cite{Sabidussi_1961} in 1961 and Hedrl\'{\i}n and Pultr \cite{Hedrlin_Pultr_1964} in 1964. Homomorphisms between undirected graphs turned out to be a powerful instrument in many fields of pure and applied mathematics of which only a few are mentioned here.

Looking at undirected graphs without loops, the pure existence of a homomorphism from $G$ to $H$ may indicate important structural features of $G$ or $H$; a wealth of topics and examples is contained in the outstanding textbook of Hell and Ne\v{s}et\v{r}il \cite{Hell_Nesetril_2004}. Many of the results even have a dual structure in the sense that for undirected graphs $G, H$, and $H'$ without loops, the existence of a homomorphism from $G$ to $H$ ist equivalent to the non-existence of a homomorphism from $H'$ to $G$. (For directed graphs, corresponding results can be found in the book of Bang-Jensen and Gutin \cite{BangJensen_Gutin_2009}.) Also homomorphisms between undirected graphs with loops are of interest in different fields of mathematics, computer science \cite{Borgs_etal_2006,Brucker_2007}, and statistical physics \cite{Widom_Rawlinson_1970,Brightwell_Winkler_1999,
Freedman_etal_2007,Borgs_etal_2008,Borgs_etal_2012} where homomorphisms describe admissible states of physical structures.

Also for posets, homomorphisms may carry deep information about structure. Birkhoff \cite{Birkhoff_1942} conjectured already in 1942, that $\H(P,R) \simeq \H(P,S)$ implies $R \simeq S$ for finite posets $P, R$, and $S$. This problem, called the ``cancellation problem'', moved into focus at the end of the seventies. Many authors contributed to its solution (a survey is contained in \cite{McKenzie_2003}), but finally, McKenzie succeeded in putting the keystone onto all the hard efforts: in 1999 and 2000, he proved \cite{McKenzie_1999,McKenzie_2000} that indeed $\H(P,R) \simeq \H(P,S)$ implies $R \simeq S$ for finite poests $P, R, S$. In 2003, McKenzie \cite{McKenzie_2003} published an additional paper about this subject, using a different approach for the proof. Another important structural topic in theoretical and applied order theory are homomorphisms with fixed points \cite{Davey_Priestley_2012,Schroeder_2016}.

Even the pure number of homomorphisms may allow conclusions to be drawn about structure. The graph parameters which can be expressed as numbers of homomorphisms into weighted graphs have been characterized by Freedman et al.\ \cite{Freedman_etal_2007} in 2007.
The still open reconstruction conjecture asks in different fields of graph theory \cite{Hell_Nesetril_2004,Schroeder_2016} if indeed two objects with at least four vertices are isomorphic if all numbers of embeddings of certain subgraphs into them are equal. A classical result is the Theorem of Lov\'{a}sz \cite{Lovasz_1967} from 1967 which states that numbers of homomorphisms distinguish non-isomorphic ``relational structures''. For our purpose, the following specifications for digraphs are relevant:

\begin{theorem}[Lov\'{a}sz \cite{Lovasz_1967}] \label{theo_Lovasz_original}
Let $\fC$ be the class of finite digraphs or the class of finite posets. Then, for $R, S \in \fC$
\begin{align*}
R & \simeq S \\
\Leftrightarrow \; \; \; \# \H(G,R) & = \# \H(G,S) \; \mytext{for every } G \in \fC.
\end{align*}
For the class of finite posets, the equivalence holds also if we replace the homomorphism sets by the sets of strict order homomorphisms.
\end{theorem}

\begin{figure} 
\begin{center}
\includegraphics[trim = 70 710 190 70, clip]{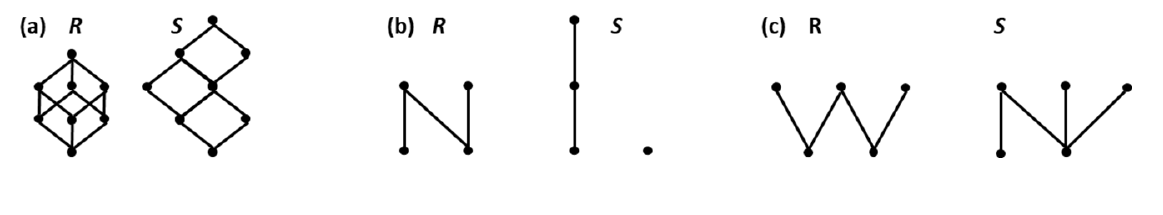}
\caption{\label{fig_Intro} The Hasse-diagrams of three pairs of posets $R$ and $S$ with $\# \H(P,R) \leq \# \H(P,S)$ for every finite poset $P$.}
\end{center}
\end{figure}

For the class of digraphs, a short and simple proof of the theorem is given in \cite{Hell_Nesetril_2004} which - with minor modification - works for posets, too.

The infinite vector $\L(H) \equiv ( \# \H(G,H) )_{G \in \fC}$ is called the {\em Lov\'{a}sz-vector of $H$}.  In the last two decades, topics related to it have found interest in connection with vertex and edge weights \cite{Borgs_etal_2006,Lovasz_2006,Freedman_etal_2007,
Borgs_etal_2008,Lovasz_Szegedy_2008,Schrijver_2009,
Borgs_etal_2012,Cai_Govorov_2020}. In the field of undirected graphs, Dvo\v{r}\'{a}k \cite{Dvorak_2010} investigated in 2010 proper sub-classes $\fU'$ of undirected graphs for which the partial Lov\'{a}sz-vector $( \# \H(G,H) )_{G \in \fU'}$ is still able to distinguish graphs; the distinguihing power of the vector $( \# \H(G,H) )_{H \in \fU'}$ is treated by Fisk \cite{Fisk_1995} in 1995.

The subject we are dealing with in this paper is the pointwise less-equal-relation between partial Lov\'{a}sz-vectors of digraphs:
\begin{quote}
{\bf Question:} Given a class $\fD'$ of digraphs, what is it in the structure of digraphs $R$ and $S$ that enforces
\begin{equation} \label{fragestellung}
\# \H(G,R) \leq \# \H(G,S) \; \mytext{\em for every } G \in \fD' \mytext{\em ?}
\end{equation}
\end{quote}

As examples, Figure \ref{fig_Intro} shows the Hasse-diagrams of three pairs of posets $R$ and $S$ with $\# \H(P,R) \leq \# \H(P,S)$ for every finite poset $P$; the proof for the pair in Figure \ref{fig_Intro}(a) in \cite{aCampo_2018} was the starting point of this investigation.

Besides of the mathematical relevance of the question, it is also of interest in application. The execution of a process in a converting facility (IT-system, logistic center, factory) can in many cases be described by a homomorphism \cite{Brucker_2007,Drozdowski_2009,
Hell_Nesetril_2004}. The converting facility is given by a reflexive digraph $R$. The vertices represent the central converting units (computers or processors, stockrooms, machine halls), and the arcs represent transport lines between them (communication lines, transportation belts, tracks and streets). Similarly, a process $G$ to be executed in a converting facility can often be described by a transitive digraph: the vertices are the process steps, and an arc from $v$ to $w$ means that $v$ has to be executed before $w$. Each homomorphism from $G$ to $R$ describes a way how the process $G$ can be executed in the facility $R$. A facility $S$ with $\# \H(G,R) \leq \# \H(G,S)$ provides thus more flexibility in execution of $G$, which may be a value in itself (e.g., in IT-systems \cite{Simon_2018}) or may enable process optimization. Moreover, the overall-structure of $S$ can be advantageous: in Figure \ref{fig_Intro}(a), $S$ has less direct lines than $R$.

\begin{figure} 
\begin{center}
\includegraphics[trim = 70 550 250 70, clip]{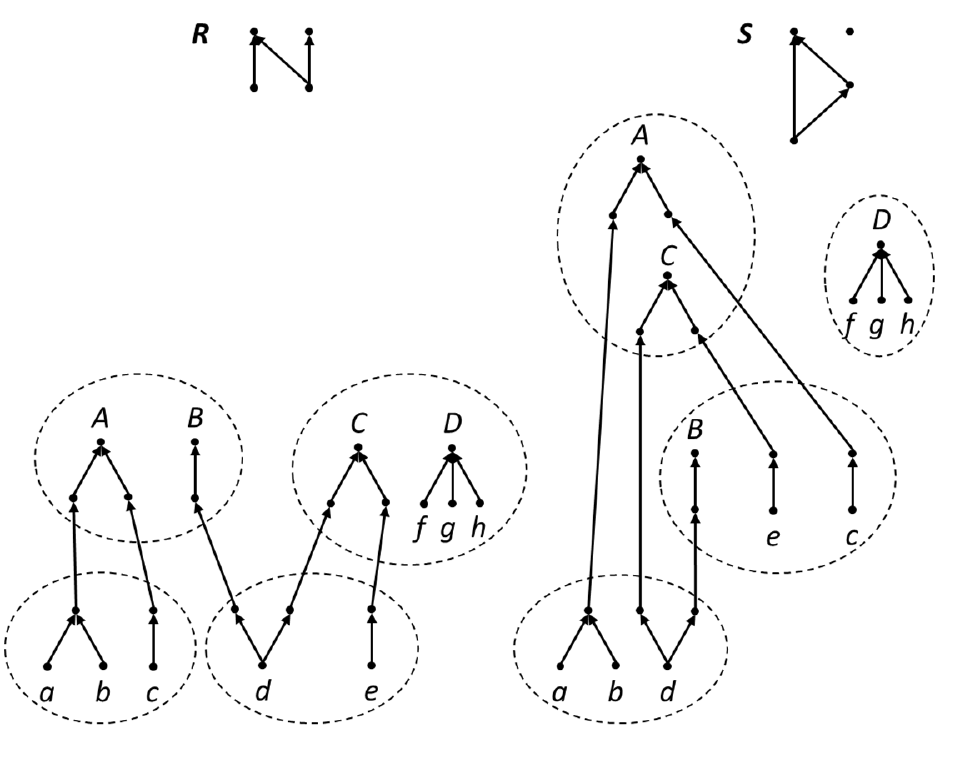}
\caption{\label{fig_ProdExample} Realization of a process $G$ on different IT-systems $R$ and $S$ (units encircled with dotted lines). Trivial arcs in the diagrams of $G, R$, and $S$ are omitted, and the sub-graph of the pre-processing of data source $d$ has been flipped horizontally in $S$ in order to avoid line crossings. Explanations in text.}
\end{center}
\end{figure}

The transfer of the execution of $G$ from $R$ to $S$ can be described by a mapping $\rho : \H(G,R) \rightarrow \H(G,S)$. However, given $\xi \in \H(G,R)$, we are in practice not free in the choice of $\rho(\xi) \in \H(G,S)$. Figure \ref{fig_ProdExample} shows on the left a - say - IT-system $R$ with N-shaped architecture on which a process $G$ is running consisting of four parallel tasks $A, B, C$, and $D$, accessing data sources $a, \ldots , h$. For the tasks $A, B$, and $C$, the homomorphism $\xi : G \rightarrow R$ in the figure sends the pre-processing of the data and the main calculation to different units of $R$. If we assume that this separation has been chosen by purpose, it should be preserved by $\rho(\xi)$. Mathematically, the pre-processings and main calculations are connectivity components of the pre-images $\myurbild{\xi}(r)$, $r \in R$, and they should be preserved by $\rho(\xi)$.

In the right part of Figure \ref{fig_ProdExample}, a possible choice for $\rho(\xi) : G \rightarrow S$ is shown, $S$ as indicated on top. Compared with $R$, the allocation of process steps to common units has heavily changed in $S$; however, every pre-processing and every main calculation has been transferred as a whole.

Our subject is thus the following: Let $R$ and $S$ be digraphs and let $\fD'$ be a class of digraphs. For every $G \in \fD'$, $\xi \in \H(G,R)$, $v \in V(G)$, let $\gxiv$ be the connectivity component of $v$ in $\myurbild{\xi}(\xiv)$. Our main question makes us interested in the existence of a one-to-one mapping $\rho_G : \H(G,R) \rightarrow \H(G,S)$ for every $G \in \fD'$, and additionally we want $\rho_G$ fulfilling $\gxiv = \grgxiv$ for every $\xi \in \H(G,R)$, $v \in G$, $G \in \fD'$. We write $R \strG S$ iff such a $\rho$ exists. Trivially, $R \strG S$ implies \eqref{fragestellung}.

The preservation of connectivity components of pre-images may look like an additional difficulty posed upon a question difficult enough in itself. However, in fact it is a regularity condition making things manageable by introducing structure. In this article, we develop several criteria for the relation $R \strG S$, hence of \eqref{fragestellung}. A study about structural properties of $R$ and $S$ enforcing $R \strG S$ is under preparation \cite{aCampo_toappear_1}.

After recalling common terms and notation in Section \ref{subsec_notation}, connectivity components of pre-images are treated in Section \ref{subsec_connectivity}. Section \ref{subsec_Def_HomSchemes} is dedicated to the definition of the relation $R \strG S$, and the first main result is derived in Section \ref{subsec_GschemeOnStrict}: We show in Theorem \ref{theo_GschemeOnStrict} that, for certain sub-classes $\fD'$ of digraphs and $R \in \fD'$, the relation  $R \strG S$ is equivalent to $\# \S(G,R) \leq \# \S(G,S)$ for all $G \in \fD'$, where $\S(G,H)$ denotes the set of {\em strict} homomorphisms from $G$ to $H$, i.e., of homomorphisms from $G$ to $H$ mapping all proper arcs of $G$ to proper arcs of $H$. A direct consequence is that for these sub-classes $\fD'$ and digraphs $R$,
\begin{align*}
\# \S(G,R) & \leq \# \S(G,S) \; \; \mytext{for all} \; G \in \fD' \\
\Rightarrow \quad \# \H(G,R) & \leq \# \H(G,S) \; \; \mytext{for all} \; G \in \fD'. 
\end{align*}
Examples for such classes $\fD'$ are the class of posets, the class of digraphs with antisymmetric transitive hull, and the full class of digraphs. Additional restrictions of $\fD'$ are possible, e.g., with respect to the maximal number of vertices or arcs. In Proposition \ref{prop_Gscheme_Hscheme}, we see that, in many cases, the relation $ R \strG S $ can be proven by showing $\# \S(G,R) \leq \# \S(G,S)$ for digraphs $G$ contained in a proper subset of $\fD'$.

The second main result is contained in Section \ref{subsec_rearrangement}. We develop a method how to rearrange a digraph $R$ fulfilling certain conditions in such a way that $R \strG S$ holds for the resulting digraph $S$ (Theorem \ref{theo_rhoxi}). We apply the method on partially ordered sets in Section \ref{subsec_posets} and construct nine pairs $R$ and $S$ of posets with $R \strG S$, including the examples in Figure \ref{fig_Intro}. Finally, in Section \ref{subsec_ugraphs}, we transfer our concepts to undirected graphs and show that parts of our results are valid also for them.

Besides of its theoretical value, Theorem \ref{theo_GschemeOnStrict} makes the investigation of the relation $\strG$ considerably simpler because it is much easier to work with strict homomorphisms than with general ones. Already this article is benefitting from it: the description of the ``exceptional'' vertex sets $U_\xi$ and $U'_\zeta$ in Definition \ref{def_Uxi} and the proof of Theorem \ref{theo_rhoxi} would become much more complicated if we had to deal with non-strict homomorphisms also.

\section{Preparation} \label{sec_preparation}

\subsection{Basics and Notation} \label{subsec_notation}

A {\em (finite) directed graph} or {\em digraph} $G$ is an ordered pair $(V(G),A(G))$ where $V(G)$ is a non-empty, finite set and $A(G) \subseteq V(G) \mytimes V(G)$ is a binary relation on $V(G)$. As usual we write $vw$ for an ordered pair $(v,w) \in V(G) \mytimes V(G)$. We call the elements of $V(G)$ the {\em vertices} of $G$ and the elements of $A(G)$ the {\em arcs} of $G$. A digraph $G$ is called {\em reflexive},  or {\em symmetric}, or {\em antisymmetric}, etc., iff the relation $A(G)$ has the respective property. A reflexive, antisymmetric, transitive digraph is called a {\em partially ordered set} or simply a {\em poset}. For a digraph $G$ and a non-empty set $X \subseteq V(G)$, the {\em digraph $G \vert_X$ induced on $X$} is defined as $( X, A(G) \cap ( X \mytimes X) )$. The {\em direct sum} $G + H$ of digraphs and the {\em ordinal sum} $P \oplus Q$ of posets are defined as usual.

Vertices $v, w \in V(G)$ are called {\em adjacent} iff $vw \in A(G)$ or $wv \in A(G)$. The {\em (open) neighborhood} $N_G(v)$ of $v \in V(G)$ is the set of all $w \in V(G) \setminus \{ v \}$ adjacent to $v$. Furthermore,
\begin{align*}
\Nin_G(v) & \; \equiv \; \mysetdescr{ w \in N_G(v) }{ wv \in A(G) }, \\
\Nout_G(v) & \; \equiv \; \mysetdescr{ w \in N_G(v) }{ vw \in A(G) }.
\end{align*}
An arc $vw \in A(G)$ with $v = w$  is called a {\em loop}; we collect all possible loops of $G$ in the {\em diagonal (relation)} $\Delta_G \equiv \mysetdescr{(v,v)}{v \in V(G)}$. $G^* \equiv (V(G), A(G) \setminus \Delta_G)$ is the digraph $G$ {\em with loops removed}. An arc $vw \in A(G)$ is called {\em proper} iff it is not a loop, i.e., iff $v \not= w$.

A sequence $W = v_0, \ldots , v_I$ of vertices of $G$ with $I \in \myN$ is called a {\em walk} iff $v_{i-1}v_i \in A(G)$ for all $1 \leq i \leq I$. The walk $W$ {\em starts in $v_0$} and {\em ends in $v_I$}. In the case of $v_0 = v_I$, the walk is {\em closed}. A walk $v_0, \ldots , v_I$ is {\em trivial} iff $v_i = v_0$ for all $1 \leq i \leq I$. A digraph is {\em acyclic} iff it does not contain a closed walk.

The term {\em convex} has different meanings in graph theory \cite[Chapter 1.4]{Pelayo_2013}. We use it as in order theory: For a digraph $G$, we call a subset $X \subseteq V(G)$ {\em convex} iff every walk starting and ending in $X$ runs totally in $X$.

Let $G$ be a digraph. With $\T$ denoting the set of all transitive relations $T \subseteq V(G) \mytimes V(G)$ with $A(G) \subseteq T$, the {\em transitive hull} $( V(G), \cap \T )$ of $G$ is the digraph with vertex set $V(G)$ and the (set-theoretically) smallest transitive arc set containing $A(G)$. If $G$ is reflexive, then its transitive hull is reflexive, too; however, antisymmetry is in general not preserved. $vw \in \cap \T $ is equivalent to the existence of a walk $z_0, \ldots , z_I $ in $G$ with $v = z_0$ and $w = z_I$. 

Given digraphs $G$ and $H$, we call a mapping $\xi : V(G) \rightarrow V(H)$ a {\em homomorphism} from $G$ to $H$ iff $\xi(v) \xi(w) \in A(H)$ for all $vw \in A(G)$. For such a mapping, we write $\xi : G \rightarrow H$, and we collect the homomorphisms in the set
\begin{align*}
\H(G,H) & \; \equiv \; \mysetdescr{ \xi : V(G) \rightarrow V(H) }{ \xi \mytext{ is a homomorphism} }.
\end{align*}

Every homomorphism $\xi : G \rightarrow H$ maps loops in $G$ to loops in $H$, but proper arcs of $G$ can be mapped to both, loops and proper arcs of $H$. We call a homomorphism from $G$ to $H$ {\em strict} iff it maps all proper arcs of $G$ to proper arcs of $H$.
\begin{align*}
\S(G,H) & \; \equiv \; \H(G,H) \cap \H(G^*,H^*)
\end{align*}
is the set of strict homomorphisms from $G$ to $H$. The set $\H(G^*,H^*) \setminus \H(G,H)$ contains all homomorphisms from $G^*$ to $H^*$ which map a vertex belonging to a loop in $G$ to a vertex of $H$ not belonging to a loop. A mapping $\xi : V(G) \rightarrow V(H)$ is thus a strict homomorphism, iff it maps loops in $G$ to loops in $H$ and proper arcs of $G$ to proper arcs of $H$. If $H$ is reflexive, then $\S(G,H) = \H(G^*,H^*)$, and for posets $P$ and $Q$, the set $\S(P,Q) = \H(P^*,Q^*)$ is the set of strict order homomorphisms from $P$ to $Q$.

We need symbols for several classes of digraphs. $\fD$ is the class of all digraphs with finite non-empty vertex set, and $\fP \subset \fD$ is the class of all finite posets.
\begin{align*}
\fP^* & \; \equiv \; \mysetdescr{ P^* }{ P \in \fP }
\end{align*}
is the class of posets with loops removed, i.e., the class of finite irreflexive antisymmetric transitive digraphs. $\fP^*$ is of interest for us, because every result about homomorphism sets $\H(P,Q)$ with $P, Q \in \fP^*$ directly translates into a result about the sets of strict order homomorphisms between posets and vice versa. For example, the addendum in Theorem \ref{theo_Lovasz_original} says that the stated equivalence also holds for $\fC = \fP^*$.

The following class will play an important role:
\begin{align*}
\fTa & \; \equiv \; \mysetdescr{ G \in \fD }{ G^* \mytext{is acyclic} }.
\end{align*}
Equivalently, $\fTa$ can be characterized as the class of digraphs with antisymmetric transitive hull (that is the reason for the choice of the symbol $\fTa$), as the class of digraphs in which every closed walk is trivial, or as the class of subgraphs of posets. In particular, all digraphs in $\fTa$ are antisymmetric, and for $G \in \fTa$, also the transitive hull of $G$ belongs to $\fTa$. Important sub-classes contained in $\fTa$ are $\fP$ and $\fP^*$ and (trivially) the class of acyclic digraphs.

For every class $\fD' \subseteq \fD$, we denote by $\fD'_r$ a representative system of $\fD'$ with respect to isomorphism. $\fD'_r$ is always a set.

Additionally, we use the following notation from set theory:
\begin{align*}
\myNk{0} & \equiv  \emptyset, \\
\myNk{n} & \equiv  \{ 1, \ldots, n \} \mytext{for every} n \in \myN.
\end{align*}

$\A(X,Y)$ is the set of mappings from $X$ to $Y$. For $f \in \A(X,Y)$ and $X' \subseteq X$, we write $f \vert_{X'}$ for the pre-restriction of $f$ to $X'$. Furthermore, we use the symbol $\myurbild{f}(Y')$ for the pre-image of $Y' \subseteq Y$ under $f$; for $y \in Y$, we simply write $\myurbild{f}(y)$ instead of $\myurbild{f}( \{y\} )$. However, in Theorem \ref{theo_rhoxi}, we use the symbol $\beta^{-1}$ also for the inverse of a bijective mapping $\beta$.

Finally, we use the {\em Cartesian product}. Let $\I$ be a non-empty set, and let $M_i$ be a non-empty set for every $i \in \I$. Then the Cartesian product of the sets $M_i, i \in \I$, is defined as
\begin{eqnarray*}
\prod_{i \in \I} M_i & \; \equiv \; & 
\mysetdescr{ f \in \A \big( \I, \bigcup_{i \in \I} M_i \big)}{ f(i) \in M_i \mytext{for all} i \in \I }.
\end{eqnarray*}

\subsection{Connectivity} \label{subsec_connectivity}

\begin{definition} \label{def_connected}
Let $G \in \fD$, $X \subseteq V(G)$, and $v, w \in X $. We say that $v$ and $w$ are {\em connected in $X$}, iff $v = w$ or, in the case of $v \not= w$, there exist $z_0, z_1, \ldots , z_I \in X$, $I \in \myN$, with $v = z_0$, $w = z_I$, and $ z_{i-1}$, $z_i$ being adjacent for all $i \in \myNk{I}$. In this case, we call $z_0, \ldots , z_I$ a {\em line connecting $v$ and $w$}. We define for all $v \in X \subseteq V(G)$
\begin{align*}
\gamma_X(v) & \equiv \mysetdescr{ w \in X }{ v \mytext{and} w \mytext{are connected in} X}.
\end{align*}
\end{definition}
The following corollary has originally been formulated for partially ordered sets in \cite{aCampo_2018}; however, it also holds for digraphs:
\begin{corollary} \label{coro_gamma_def}
Let $G \in \fD$ and $X \subseteq V(G)$. The relation ``connected in $X$'' is an equivalence relation on $X$ with partition $\mysetdescr{\gamma_X(x)}{x \in X}$. For $x \in X \subseteq X' \subseteq V(G)$ we have
\begin{align}
\gamma_X(x) & \subseteq \gamma_{X'}(x), \label{gamma_ABAB} \\
\gamma_X(x) & = \gamma_{\gamma_X(x)}(x). \label{gamma_ggc}
\end{align}
\end{corollary}

The sets $\gamma_{V(G)}(v), v \in V(G)$, are called the {\em connectivity components} of $G$. Every digraph is the direct sum of its connectivity components. A digraph $G \in \fD$ is {\em connected} iff $\gamma_{V(G)}(v) = V(G)$ for a $v \in V(G)$ (the choice of $v$ is arbitrary). A subset $X \subseteq V(G)$ is called {\em connected (in $G$)} iff the digraph $G \vert_X$ induced on $X$ is connected.

The following definition is one of the central ones in this paper:
\begin{definition}
Let $G \in \fD$, let $X$ be a set, and let $\xi \in \A(V(G),X)$ be a mapping. We define for all $v \in V(G)$
\begin{align*}
\gxiv & \equiv \gamma_{\myurbild{\xi}( \xiv )}( v ).
\end{align*}
\end{definition}
$\gxiv$ contains thus $v$ and all vertices, $v$ is connected with in $\myurbild{\xi}( \xiv )$. According to \eqref{gamma_ggc}, $\gxiv = \gamma_{\gxiv}(v)$ for every $v \in V(G)$: for $w \in \gxiv$, $w \not= v$, there is a line connecting $v$ and $w$ in $\myurbild{\xi}(\xiv)$ that runs totally in $\gxiv$.

\begin{corollary} \label{coro_gxix_strict}
Let $G, H \in \fD$. Then for every homomorphism $\xi \in \H(G,H)$
\begin{equation*} 
\xi \; \mytext{strict} \quad \Leftrightarrow \quad \forall \; v \in V(G) \; : \; \gxiv = \{ v \}.
\end{equation*}
\end{corollary}
\BP If $\xi$ is strict, then $\xiw \not= \xiv$ for every $w \in N_G(v) $, thus $\gxiv = \{ v \}$. On the other hand, let $v, w \in V(G^*) = V(G)$ with $vw \in A(G^*)$. $\gxiv = \{ v \}$ yields $\xiw \not= \xiv$, hence $\xiv \xiw \in A(H^*)$.

\EP

\begin{lemma} \label{lemma_Gxix_subsetnoteq_Gzetax}
Let $G, H, H' \in \fD$, let $\xi \in \H(G,H)$ and $\zeta \in \H(G,H')$ be homomorphisms, and let $\gxiv \subseteq \gzev$ for a vertex $v \in V(G)$. Then $\gxiv \subset \gzev$ iff there are $a, b \in \gzev$ for which $ab$ is a proper edge in $G$ and $\xi(a) \xi(b)$ is a proper edge in $H$.
\end{lemma}
\BP Let $w \in \gzev \setminus \gxiv$. The points $w$ and $v$ are connected by a line in $\gzev$, and on this line there exist adjacent vertices $a$ and $b$ with $a \in \gxiv$, $b \in \gzev \setminus \gxiv$. We have $a = v$, or $a$ and $v$ are connected by a line in $\gxiv$. In both cases, $b$ and $v$ are connected by a line in $\{ b \} \cup \gxiv$, and $b \notin \gxiv$ means $\xi(b) \not= \xiv = \xi(a)$, and $\xi(a) \xi(b)$ or $\xi(b) \xi(a)$ is a proper arc of $H$, depending on $ab \in A(G)$ or $ba \in A(G)$.

On the other hand, let $a, b \in \gzev$ with $a b$ and $\xi(a) \xi(b)$ being proper arcs. Due to $\gxiv \subseteq \myurbild{\xi}(\xiv)$, we conclude that at least one of the vertices $a, b$ does not belong to $\gxiv$.

\EP

\begin{corollary} \label{coro_G_sigma_xi}
Let $G, H, H'$ be finite digraphs, and let $\xi \in \H(G,H)$ and $\sigma \in \H(H,H')$. Then $\gxiv \subseteq \mygxiv{\sigma \circ \xi}{v}$ for all $v \in V(G)$. Equality holds for all $v \in V(G)$, if $\sigma \vert_{\xi[V(G)]}$ is strict.
\end{corollary}
\BP We have $\myurbild{ \xi }( \xiv) \subseteq \myurbild{ \xi }(  \myurbild{ \sigma } ( \sigma( \xiv)) = \myurbild{ { ( \sigma \circ \xi ) } }(( \sigma \circ \xi)(v))$ for every $v \in V(G)$, and $\gxiv \subseteq \mygxiv{\sigma \circ \xi}{v}$ follows with \eqref{gamma_ABAB}.

Assume $\gxiv \subset \mygxiv{\sigma \circ \xi}{v} $ for a $v \in V(G)$. Lemma \ref{lemma_Gxix_subsetnoteq_Gzetax} delivers the existence of adjacent $a, b \in \mygxiv{\sigma \circ \xi}{v}$ for which $\xi(a) \xi(b)$ is a proper edge in $H$. $a, b \in \mygxiv{\sigma \circ \xi}{v}$ implies $\sigma( \xi(a) ) = \sigma( \xi(b) )$, and $\sigma$ is not strict on $\xi[V(G)]$.

\EP

\section{Strong $\bf \Gamma$-schemes} \label{sec_strongHomSchemes}

\subsection{Definition of strong $\bf \Gamma$-schemes} \label{subsec_Def_HomSchemes}

Let $R, S \in \fD$ and $\fD' \subseteq \fD$. Assume that there exists a one-to-one homomorphism $\sigma : R \rightarrow S$. Then, for every $G \in \fD'$ with $\H(G,R) \not= \emptyset$, we get a one-to-one mapping $r_G : \H(G,R) \rightarrow \H(G,S)$ by setting for every $\xi \in \H(G,R)$
\begin{align} \label{eq_rho_sigma}
r_G(\xi) & \equiv \sigma \circ \xi.
\end{align}
A one-to-one homomorphisms from $R$ to $S$ delivers thus a ``natural'' (or: trivial) example for $ \# \H(G,R) \leq \# \H(G,S)$ for every $G \in \fD'$. Furthermore, according to Corollary \ref{coro_G_sigma_xi}, we have
\begin{align*}
\mygxiv{r_G(\xi)}{v} & = \gxiv
\end{align*}
for every $G \in \fD', \xi \in \H(G,R), v \in V(G)$. These properties of \eqref{eq_rho_sigma} are taken up in the following definition; we have to use the Cartesian product in it in order to make the letter ``$\rho$'' being a meaningful mathematical object.

\begin{definition} \label{def_RS_scheme}
Let $\fD' \subseteq \fD$ be a sub-class of digraphs. For $R, S \in \fD$, we call a mapping
\begin{align*}
\rho & \; \in \prod_{G \in \fD'_r} \A( \H(G,R), \H(G,S) )
\end{align*}
a {\em Hom-scheme from $R$ to $S$ (with respect to $\fD'$)}, and we call it {\em strong} iff $\rho_G : \H(G,R) \rightarrow \H(G,S)$ is one-to-one for every $G \in \fD'_r$. We say that a Hom-scheme $\rho$ from $R$ to $S$ is a {\em $\Gamma$-scheme}, iff
\begin{align} \label{def_grxix_gxix}
\grgxixv & = \gxiv
\end{align}
for every $G \in \fD'_r, \xi \in \H(G,R), v \in V(G)$. We write $R \strG S$ iff a strong $\Gamma$-scheme from $R$ to $S$ exists. If $G$ is fixed, we write $\rho(\xi)$ instead of $\rho_G(\xi)$.
\end{definition}

Here as in the following, it does not matter if there is a $G \in \fD'$ with $\H(G,D) = \emptyset$; in this case, $\rho_G = ( \emptyset, \emptyset, \H(G,S))$. The Hom-scheme $r$ in \eqref{eq_rho_sigma} induced by a one-to-one homomorphism $\sigma$ is always a strong $\Gamma$-scheme. If the mapping $\sigma$ in \eqref{eq_rho_sigma} is only strict, then $r$ is a still a $\Gamma$-scheme according to Corollary \ref{coro_G_sigma_xi}, but not necessary a strong one.

There exists a strong Hom-scheme $\rho$ from $R$ to $S$ iff
\begin{equation*}
\# \H(G,R) \; \leq \; \# \H(G,S) \; \; \mytext{for all} \; G \in \fD',
\end{equation*}
and $\rho$ is a strong $\Gamma$-scheme iff it additionally obeys the regularity condition \eqref{def_grxix_gxix} in mapping $\H(G,R)$ to $\H(G,S)$ for every $G \in \fD'_r$. This regularity condition is plausible if we regard a Hom-scheme as a technical apparatus which assigns to every $\xi \in \H(G,R)$ a well-fitting $\rho(\xi) \in \H(G,S)$. If we allow $\gxiv \subset \grxiv$ for $v \in V(G)$, then Lemma \ref{lemma_Gxix_subsetnoteq_Gzetax} tells us that $\rho(\xi)$ preserves the structure of $G$ around $v$ worse than $\xi$, which is not satisfying. And in the case $\gxiv \not\subseteq \grxiv$, $\rho(\xi)$ has to re-distribute the points of $\gxiv \setminus \grxiv \subseteq \gxiv \setminus \{ v \}$ in $S$. Because the sets $\gxiv \setminus \grxiv$ can be arbitrarily complicated, this re-distribution process may require many single case decisions, which is out of the scope of a technical apparatus.

The relation ``$\strG$'' between digraphs is always reflexive and transitive. For $\fD'$ being one of the classes $\fD, \fP$, or $\fP^*$, it even defines a partial order on $\fD'_r$: For $R, S \in \fD'_r$ with $R \strG S$ and $S \strG R$, we have $\# \H(G,R) = \# \H(G,S)$ for all $G \in \fD'_r$ which is equivalent to $R \simeq S$ according to Theorem \ref{theo_Lovasz_original} (for $\fP^*$, use the addendum).

\subsection{Strong $\bf \Gamma$-schemes and strict homomorphisms} \label{subsec_GschemeOnStrict}

In this section, we prove

\begin{theorem} \label{theo_GschemeOnStrict}
Let $R \in \fD$ and $\fD' \subseteq \fD$. Then, for all $S \in \fD$, the equivalence
\begin{align} \label{R_homG_S} 
R & \strG S \; \; \mytext{with respect to } \fD'\\ \label{SPR_leq_SPS} 
\Leftrightarrow \quad \# \S(G,R) & \leq \# \S(G,S) \; \; \mytext{for all} \; G \in \fD',
\end{align}
and the implication
\begin{align*}
\# \S(G,R) & \leq \# \S(G,S) \; \; \mytext{for all} \; G \in \fD' \\
\Rightarrow \quad \# \H(G,R) & \leq \# \H(G,S) \; \; \mytext{for all} \; G \in \fD'. 
\end{align*}
hold if
\begin{align*}
R & \in \fD' = \fD, \\
R & \in \fTa \subseteq \fD' \subseteq \fD, \\
\mytext{or} \quad R & \in \fD' \mytext{ with } \fD' = \fP \mytext{ or } \fD' = \fP^*.
\end{align*}
\end{theorem}
The implication $\eqref{R_homG_S} \Rightarrow \eqref{SPR_leq_SPS}$ does not depend at all on the choice of $R$ or $\fD'$, and its proof is simple. As stated in Corollary \ref{coro_gxix_strict}, a homomorphism $\xi \in \H(G,H)$ is strict iff $\gxiv = \{ v \}$ for every $v \in V(G)$. For a $\Gamma$-scheme $\rho$ from $R$ to $S$, we have $\grxiv = \gxiv$ for every $\ariestr$, thus $\rho_G[ \S(G,R) ] \subseteq \S(G,S)$ for every $G \in \fD'_r$, and \eqref{SPR_leq_SPS} follows if $\rho$ is strong.

It is thus the direction $\eqref{SPR_leq_SPS} \Rightarrow \eqref{R_homG_S}$ (which also yields the implication) which is of interest in the theorem. Due to $\fP, \fP^* \subset \fTa \subset \fD$, the choice of $R$ and $\fD'$ becomes more and more specialized in the three cases. The theorem states that in all three cases, the class $\fD'$ is so large and the structure of $R$ is so rich that - for all digraphs $S$ - the relation $ \# \S(G,R) \leq \# \S(G,S) $ for all $G \in \fD'$ enforces $R \strG S$ and hence $ \# \H(G,R) \leq \# \H(G,S)$ for all $G \in \fD'$.

For the proof of $\eqref{SPR_leq_SPS} \Rightarrow \eqref{R_homG_S}$, we modify the well-known mathematical approach to let a mapping $f : X \rightarrow Y$ factorize over the set $X / f \equiv \mysetdescr{ \myurbild{f}(f(x)) }{ x \in X}$ of its pre-images: $f = \iota_f \circ \pi_f$ with canonical mappings $\pi_f : X \rightarrow X / f$, $\iota_f : X / f \rightarrow Y$. In this section, we let a homomorphism factorize over a {\em refinement} of the set of its pre-images, a refinement consisting of the {\em connected} sets $\gxiv$:

\begin{definition} \label{def_pixi_ioxi}
Let $G, H \in \fD$. We define for every $ \xi \in \H(G,H)$ the digraph $\Gxi$ by
\begin{align*}
V(\Gxi) & \; \equiv \; \mysetdescr{ \gxiv }{ v \in V(G) }, \\
A(\Gxi) & \; \equiv \; \mysetdescr{ ( \fa, \fb ) \in V(\Gxi) \mytimes V(\Gxi) }{
\exists \; a \in \fa, b \in \fb \; : \; a b \in A(G)}.
\end{align*}
Additionally, we define the mappings
\begin{align*}
\pi_\xi : G & \rightarrow \Gxi \\
v & \mapsto \gxiv, \\
\iota_\xi : \Gxi & \rightarrow H \\
\gxiv & \mapsto \xiv.
\end{align*}
\end{definition}

$V(\Gxi)$ is a partition of $V(G)$ consisting of connected sets, and $\iota_\xi$ is a well-defined mapping because $\xi$ is constant on every $\fc \in \Gxi$. Obviously, $\pi_\xi \in \H(G, \Gxi)$ and $\iota_\xi \in \H(\Gxi, H)$ with $\xi = \iota_\xi \circ \pi_\xi$ for all $\xi \in \H(G,H)$. Moreover:

\begin{corollary} \label{coro_iota_strict}
$\iota_\xi$ is strict.
\end{corollary}
\BP Let $\fa, \fb \in V( \Gxi )$ with $\fa \fb \in A( \Gxi )$ and $\iota_\xi(\fa) = \iota_\xi(\fb)$. There exists $a \in \fa$, $b \in \fb$ with $ab \in A(G)$. Because $\fa$ and $\fb$ are both connected, also $\fa \cup \fb$ is connected, and $\iota_\xi(\fa) = \iota_\xi(\fb)$ yields $\xi(a) = \xi(b)$. $\xi$ is thus constant on the connected set $\fa \cup \fb$, hence $\fa \cup \fb \subseteq \mygxiv{\xi}{a} = \fa$ and $\fa \cup \fb \subseteq \mygxiv{\xi}{b} = \fb$. Therefore $\fa = \fb$, and $\iota_\xi$ is strict.

\EP

\begin{lemma} \label{lemma_Gscheme_Gamma}
For every $G, H \in \fD$, $ \xi \in \H(G,H)$, we define for every $H' \in \fD$
\begin{align*}
\Theta_{G,H'}(\xi) & \equiv \mysetdescr{ \zeta \in \H(G,H') }{ \Gzet = \Gxi }.
\end{align*}
Then, for $R, S \in \fD$, the relation $R \strG S$ with respect to $\fD' \subseteq \fD$ is equivalent to
\begin{align*}
\# \Theta_{G,R}(\xi) \leq \# \Theta_{G,S}(\xi) \; \; \mytext{for all} \; G \in \fD', \xi \in \H(G,R).
\end{align*}
\end{lemma}
\BP Let $R, S \in \fD$, $G \in \fD'$, $\xi \in \H(G,R)$, $\zeta \in \H(G,S)$. Because $V(\Gxi)$ and $V(\Gzet)$ are both partitions of $V(G)$, there is for every $v \in V(G)$ a unique $\fa \in V(\Gxi)$ and a unique $\fb \in V(\Gzet)$ with $v \in \fa = \gxiv$ and $v \in \fb = \gzev$. $\Gzet = \Gxi$ is thus equivalent to $\gxiv = \gzev$ for all $v \in V(G)$. A $\Gamma$-scheme $\rho$ from $R$ to $S$ maps thus $\Theta_{G,R}(\xi)$ to $\Theta_{G,S}(\xi)$ for every $G \in \fD'_r$, $\xi \in \H(G,R)$, and if it is strong, it does so one-to-one. On the other hand, $\# \Theta_{G,R}(\xi) \leq \# \Theta_{G,S}(\xi)$ for every $G \in \fD', \xi \in \H(G,R)$, gives trivially raise to a strong $\Gamma$-scheme, because for every $G \in \fD'_r$ with $\H(G,R) \not= \emptyset$, the set $\mysetdescr{ \Theta_{G,R}(\xi) }{ \xi \in \H(G,R) }$ is a partition of $\H(G,R)$, and the set $\mysetdescr{ \Theta_{G,S}(\xi) }{ \xi \in \H(G,R) }$ is a collection of disjoint subsets of $\H(G,S)$.

\EP

\begin{corollary} \label{coro_Gamma_pi_iota}
Let $G, H \in \fD$ and $\xi \in \H(G,H)$. Then for every $H' \in \fD$
\begin{align} \label{eq_pi_zeta}
\pi_\zeta & = \pi_\xi \quad \mytext{for all} \; \zeta \in \Theta_{G,H'}(\xi),
\end{align}
and
\begin{align} \label{equiv_pi_zeta}
\zeta_1 \not= \zeta_2 \; & \Leftrightarrow \; \iota_{\zeta_1} \not= \iota_{\zeta_2} \quad \mytext{for all} \; \zeta_1, \zeta_2 \in \Theta_{G,H'}(\xi).
\end{align}
\end{corollary}
\BP Let $\zeta \in \Theta_{G,H'}(\xi)$. Due to $\Gxi = \Gzet$, $\pi_\xi$ and $\pi_\zeta$ are both elements of $\H(G, \Gxi)$. Let $v \in V(G)$. Because $\Gxi = \Gzet$ is a partition of $V(G)$, there is a unique $\fa \in \Gxi = \Gzet$ with $v \in \fa$. We conclude $\pi_\xiv = \gxiv = \fa = \gzev = \pi_\zeta(v)$, and \eqref{eq_pi_zeta} is shown, because $v \in V(G)$ was arbitrary.

For $\zeta_1, \zeta_2 \in \Theta_{G,H'}(\xi)$, we have $\iota_{\zeta_1}, \iota_{\zeta_2} \in \H( \Gxi, H' )$. The conclusion $\zeta_1 = \zeta_2 \Rightarrow \iota_{\zeta_1} = \iota_{\zeta_2}$ is trivial. Let $\zeta_1 \not= \zeta_2$. Then $\iota_{\zeta_1} \circ \pi_{\zeta_1} = \zeta_1 \not= \zeta_2 = \iota_{\zeta_2} \circ \pi_{\zeta_2}$, and \eqref{eq_pi_zeta} delivers $\iota_{\zeta_1} \not= \iota_{\zeta_2}$

\EP

\begin{lemma} \label{lemma_Gamma_SPGxiQ}
For every $G, H \in \fD$ and $\xi \in \H(G,H)$, we have $\# \Theta_{G,H'}(\xi) = \# \S( \Gxi, H' )$ for every $H' \in \fD$.
\end{lemma}
\BP Let $\J \equiv \mysetdescr{ \iota_\zeta }{ \zeta \in \Theta_{G,H'}(\xi) }$. Corollary \ref{coro_iota_strict} yields $\J \subseteq \S(\Gxi,H')$, and \eqref{equiv_pi_zeta} delivers $ \# \Theta_{G,H'}(\xi) = \# \J$. We conclude $\# \Theta_{G,H'}(\xi) \leq \# \S( \Gxi, H' )$.

Let now $\sigma_1, \sigma_2 \in \S( \Gxi, H' )$ with $\sigma_1 \not= \sigma_2$. With $\zeta_1 \equiv \sigma_1 \circ \pi_\xi$, $\zeta_2 \equiv \sigma_2 \circ \pi_\xi$ we have $\zeta_1, \zeta_2 \in \H(G,H')$.

Let $i \in \myNk{2}$ be fixed. We want to show $\zeta_i \in \Theta_{G,H'}(\xi)$. For $v \in V(G)$, the set $\gxiv$ contains $v$ and is connected in $G$, and the mapping $\zeta_i$ is constant on $\pi_\xiv = \gxiv$. Therefore, $\gxiv \subseteq \mygxiv{\zeta_i}{v}$. In the case of ``$\subset$'', Lemma \ref{lemma_Gxix_subsetnoteq_Gzetax} delivers
$a, b \in \mygxiv{\zeta_i}{x}$ for which $a b$ and $\xi(a) \xi(b)$ are proper arcs in $G$ and $H$, respectively. But $a b \in A(G)$ means $\mygxiv{\xi}{a} \mygxiv{\xi}{b} = \pi_\xi(a) \pi_\xi(b) \in A(\Gxi)$, and $\xi(a) \not= \xi(b)$ means $\mygxiv{\xi}{a} \not= \mygxiv{\xi}{b}$. $\mygxiv{\xi}{a} \mygxiv{\xi}{b}$ is thus a proper arc in $\Gxi$, and because $\sigma_i$ is strict, we have
\begin{displaymath}
\zeta_i(a) = \sigma_i ( \pi_\xi( a ) ) = \sigma_i( \mygxiv{\xi}{a} ) \not= \sigma_i( \mygxiv{\xi}{b} ) = \sigma_i ( \pi_\xi( b ) ) = \zeta_i(b)
\end{displaymath}
in contradiction to $a, b \in \mygxiv{\zeta_i}{v}$. Therefore, $\gxiv = \mygxiv{\zeta_i}{v}$, thus $\zeta_i \in \Theta_{G,H'}(\xi)$, because $v \in V(G)$ was arbitrary.

We have $\iota_{\zeta_1} = \sigma_1$ and $\iota_{\zeta_2} = \sigma_2$. Equivalence \eqref{equiv_pi_zeta} yields $\zeta_1 \not= \zeta_2$, and $\# \S( \Gxi, H' ) \leq \# \Theta_{G;H'}(\xi)$ is shown.

\EP

The following lemma is the key for the proof of implication \eqref{SPR_leq_SPS} $\Rightarrow$ \eqref{R_homG_S} in all three cases:

\begin{lemma} \label{lemma_GDR}
Let $R, S \in \fD$, $\fD' \subseteq \fD$, and define
\begin{equation*}
\fG_{\fD'}(R) \quad \equiv \quad \mysetdescr{ \Gxi }{ \xi \in \H(G,R), G \in \fD' }.
\end{equation*}
Then
\begin{align} \label{impl_Dstrich_S}
\# \S(G,R) & \leq \# \S(G,S) \; \; \mytext{for all} \; G \in \fG_{\fD'}(R) \\ \label{impl_Dstrich_H}
\mytext{implies} \quad \quad \quad \quad \quad \quad R &\strG S \; \; \mytext{with respect to} \; \fD'.
\end{align}
\end{lemma}
\BP Let $G \in \fD'$, $\xi \in \H(G,R)$. Applying Lemma \ref{lemma_Gamma_SPGxiQ} with $H = H' = R$ yields $\# \Theta_{G,R}(\xi) = \# \S( \Gxi, R )$, and using Lemma \ref{lemma_Gamma_SPGxiQ} with $H = R, H' = S$ results in $\# \Theta_{G,S}(\xi) = \# \S( \Gxi, S )$. Assumption \eqref{impl_Dstrich_S} delivers $\# \Theta_{G,R}(\xi) \leq \# \Theta_{G,S}(\xi)$, and Lemma \ref{lemma_Gscheme_Gamma} yields \eqref{impl_Dstrich_H}.

\EP 

For the case $R \in \fD' = \fD$, the implication \eqref{SPR_leq_SPS} $\Rightarrow$ \eqref{R_homG_S} in Theorem \ref{theo_GschemeOnStrict} is a direct consequence of this lemma. For the remaining cases, we need that $\iota_\xi$ is ``strict across walks'' for $H \in \fTa$:

\begin{lemma} \label{lemma_iota}
Let $G \in \fD$, $H \in \fTa$ and $ \xi \in \H(G,H)$. If $W = \fc_0, \ldots , \fc_I$ is a walk in $\Gxi$ with $\iota_\xi( \fc_0 ) = \iota_\xi( \fc_I )$, then $W$ is a trivial walk. In particular, $\Gxi \in \fTa$.
\end{lemma}
\BP Let $\fc_0, \ldots , \fc_I$ be a walk in $\Gxi$. $\fc_{i-1} \fc_i \in A(\Gxi)$ for every $i \in \myNk{I}$ is equivalent to the existence of $v_0^+ \in \fc_0$, $v_i^-, v_i^+ \in \fc_i$ for every $i \in \myNk{{I-1}}$, and $v_I^- \in \fc_I$ with $v_{i-1}^+ v_i^- \in A(G)$ for every $i \in \myNk{I}$. Because the sets $\fc_i$ are all connected in $G$, also the set $C \equiv \bigcup_{i=0}^I \fc_i$ is connected in $G$.

We have $\xi( v_{i-1}^+ ) \xi( v_i^- ) \in A(H)$ for every $i \in \myNk{I}$, and because $\xi$ is constant on every set $\fc_i$, we have $\xi( v_i^- ) = \xi( v_i^+ )$ for every $i \in \myNk{{I-1}}$. Therefore, $W' \equiv \xi(v_0^+), \xi(v_1^+), \ldots , \xi(v_{I-1}^+), \xi(v_I^-)$ is a walk in $H$.

In the case of $\iota_\xi( \fc_0 ) = \iota_\xi( \fc_I )$, $W'$ is a closed walk in $H$, hence a trivial walk. But $\xi(v_0^+) = \xi(v_1^+) = \ldots = \xi(v_{I-1}^+) = \xi(v_I^-)$ means that $\xi$ is constant on the connected set $C$. We conclude $\fc_i = \mygxiv{\xi}{v_i^+} = \gamma_{\myurbild{\xi}(\xi(v_i^+)}(v_i^+) \supseteq C$ for every $i \in \myNk{{I-1}} \cup \{ 0 \}$, and similarly $\fc_I = \mygxiv{\xi}{v_I^-} \supseteq C$, which yields $\fc_i = C$ for every $i \in \myNk{{I}} \cup \{ 0 \}$.

$\Gxi \in \fTa$ follows a forteriori.

\EP

Now we can prove \eqref{SPR_leq_SPS} $\Rightarrow$ \eqref{R_homG_S} also for the case $R \in \fTa \subseteq \fD' \subseteq \fD$. We have $\Gxi \in \fTa$ according to the addendum in Lemma \ref{lemma_iota}, and assumption \eqref{SPR_leq_SPS} yields $ \# \S( \Gxi, R ) \leq \# \S( \Gxi, S )$ for all $\xi \in \H(G,R)$, $G \in \fD'$. Now apply Lemma \ref{lemma_GDR} again.

Transitivity of $H$ is in general not inherited by $\Gxi$, but we have
\begin{proposition} \label{prop_Txi_poset}
Let $G, H \in \fD$, $\xi \in \H(G,H)$, and let $\Txi$ denote the transitive hull of $\Gxi$. If $H$ is transitive, then
\begin{align} \label{HTxi_HGxi}
\H( \Txi, H ) & \; = \; \H( \Gxi, H ),
\end{align}
and if $H$ is additionally antisymmetric, we even have
\begin{align} \label{STxi_SGxi}
\S( \Txi, H ) & \; = \; \S( \Gxi, H ).
\end{align}
Furthermore, if $G, H \in \fP'$ with $\fP' = \fP$ or $\fP' = \fP^*$, then $\Txi \in \fP'$.
\end{proposition}
\BP `$\subseteq$'' holds in \eqref{HTxi_HGxi} and \eqref{STxi_SGxi} because of $V( \Gxi ) = V( \Txi )$ and $A( \Gxi ) \subseteq A( \Txi )$. For $\fa \fb \in A( \Txi )$, there exists a walk $W = \fc_0, \ldots , \fc_I$ in $\Gxi$ with $\fa = \fc_0$ and $\fb = \fc_I$.

Let $H$ be transitive and $\sigma \in \H(\G(\xi),H)$. The sequence $\sigma( \fc_0 ), \ldots , \sigma( \fc_I )$ is a walk in $H$, hence $\sigma( \fa ) \sigma( \fb ) \in A( H )$, and \eqref{HTxi_HGxi} is proven.

Now assume that $H$ is additionally antisymmetric (then $H \in \fTa$) and that $\sigma$ is strict. Due to Corollary \ref{coro_gxix_strict}, we have
\begin{equation*}
V( \myGxi{\sigma} ) \quad = \quad \mysetdescr{ \{ \fc \} }{ \fc \in V(\Gxi) },
\end{equation*}
thus $\myGxi{\sigma} \simeq \Gxi$. In the case of $\sigma( \fa ) = \sigma( \fb )$, the main statement in Lemma \ref{lemma_iota} applied on $\iota_\sigma : \myGxi{\sigma} \rightarrow H$ and the walk $ \{ \fc_0 \}, \ldots , \{ \fc_I \}$ in $\myGxi{\sigma}$ yields $\{ \fa \} = \{ \fb \}$, and $\sigma \in \S(\Txi, H)$ is shown.

Due to $\fP' \subset \fTa$ and Lemma \ref{lemma_iota}, we have $\G(\xi) \in \fTa$ for both choices of $\fP'$, and $\Txi$ is antisymmetric. If $G$ is reflexive, then $\Gxi$ and $\Txi$ are reflexive, too, hence $\Txi \in \fP$.

Assume $H \in \fP^*$. Due to $\iota_\xi \in \S( \Gxi, H )$, we have $\H( \Gxi, H ) \not= \emptyset$, and $\Gxi$ must be irreflexive. But an irreflexive element of $\fTa$ cannot contain any closed walk (even not a trivial one), and $\Txi$ must be irreflexive, thus an element of $\fP^*$.

\EP

For $R \in \fP'$ with $\fP' = \fP$ or $\fP' = \fP^*$, Proposition \ref{prop_Txi_poset} and  assumption \eqref{SPR_leq_SPS} yield $ \# \S( \Gxi, R ) = \# \S( \Txi, R )  \leq  \# \S( \Txi, S )$ for all $\xi \in \H(G,R), G \in \fD'$, and Lemma \ref{lemma_GDR} delivers \eqref{R_homG_S} because $\S( \Txi, S ) \subseteq \S( \Gxi, S )$ is trivial.

There are many sub-classes $\fD' \subseteq \fD$ with $\fG_{\fD'}(R) \subseteq \fD'$ for all $R \in \fD'$, e.g., the digraphs with at most $k$ vertices or at most $k$ edges. For such a sub-class $\fD'$, we can thus establish $R \strG S$ by showing \eqref{impl_Dstrich_S} for all $G \in \fD'$. Because Lemma \ref{lemma_iota} delivers $\Gxi \in \fTa$ for $R \in \fTa$, we can also choose a suitable sub-class $\fD' \subseteq \fTa$ with $\fG_{\fD'}(R) \subseteq \fD'$ for $R \in \fTa$, e.g., the class of all digraphs in $\fTa$ in which the maximal length of a walk without loops is at most $k$. For the corresponding sub-classes of $\fP$ and $\fP^*$, additionally Proposition \ref{prop_Txi_poset} has to be used.

In the following proposition we show that a strong $\Gamma$-scheme $\rho$ can, in many cases, be constructed by means of a simpler object $\sigma$ referring to sets of strict homomorphisms only:

\begin{proposition} \label{prop_Gscheme_Hscheme}
Let $R \in \fD' \subseteq \fD$ with $\fG_{\fD'}(R) \subseteq \fD'$. Then, for all $S \in \fD$, the relation $R \strG S$ is equivalent to the existence of a one-to-one mapping $\sigma_\fg : \S(\fg,R) \rightarrow \S(\fg,S)$ for every $\fg \in \fG_{\fD'}(R)$. In this case, we get a strong $\Gamma$-scheme $\rho$ from $R$ to $S$ by defining
\begin{equation} \label{rho_eta_pi}
\rho_G(\xi) \quad \equiv \quad \sigma_{\Gxi}( \iota_\xi ) \circ \pi_\xi
\end{equation}
for every $\ariegxistr$. The corresponding results hold also for $\fD' \subseteq \fP$ and $\fD' \subseteq \fP^*$, if we replace the elements of $\fG_{\fD'}(R)$ by their transitive hulls and $\Gxi$ in \eqref{rho_eta_pi} by $\Txi$.
\end{proposition}
\BP Let $R \in \fD' \subseteq \fD$ with $\fG_{\fD'}(R) \subseteq \fD'$, and let $S \in \fD$. For $R \sqsubseteq_\Gamma S$, we see $\# \S(\fg,R) \leq \# \S(\fg,S)$ for every $\fg \in \fG_{\fD'}(R)$ just as in the proof of \eqref{R_homG_S} $\Rightarrow$ \eqref{SPR_leq_SPS}.

Now assume that there exists a one-to-one mapping $\sigma_\fg : \S(\fg,R) \rightarrow \S(\fg,S)$ for every $\fg \in \fG$. According to Corollary \ref{coro_iota_strict}, we have $\iota_\xi \in \S( \Gxi, R )$, and $\rho$ in \eqref{rho_eta_pi} is well-defined. Corollary \ref{coro_G_sigma_xi} delivers for all $\ariestr$
\begin{equation*}
\mygxiv{\rho_G(\xi)}{v} \; = \; \mygxiv{ \pi_\xi }{v} \; = \; \pi_\xi(v) \; = \; \gxiv,
\end{equation*}
and $\rho$ is a $\Gamma$-scheme, thus $\rho_G[ \Theta_{G,R}(\xi) ] \subseteq \Theta_{G,S}(\xi)$.

Let $\ariegxistr$. For every $\theta, \theta' \in \Theta_{G,R}(\xi)$, we have $\pi_\theta = \pi_{\theta'}$ according to Corollary \ref{coro_Gamma_pi_iota}. $\pi_\theta$ is onto and $\sigma_{\G(\xi)}$ is one-to-one. Therefore, $\rho_G(\theta) = \rho_G(\theta')$ implies $\iota_\theta = \iota_{\theta'}$, hence $\theta = \theta'$ according to \eqref{equiv_pi_zeta} (with $H = H' = R$). $\rho_G$ is thus one-to-one on $\Theta_{G,R}(\xi)$, hence $\# \Theta_{G,R}(\xi) \leq \# \Theta_{G,S}(\xi)$, and just as at the end of the proof of Lemma \ref{lemma_Gscheme_Gamma}, we conclude that $\rho$ is a strong $\Gamma$-scheme from $R$ to $S$.

For $\fD' \subseteq \fP$ and $\fD' \subseteq \fP^*$, use Proposition \ref{prop_Txi_poset} and observe $\pi_\xi \in \H(G, \Txi)$.

\EP

Finally, we show that $\strG$ is compatible with the direct sum of digraphs; compatibility with other operations and cancellation rules are the subject of a separate paper \cite{aCampo_toappear_2}:

\begin{corollary} \label{calc_dirsum}
Let $\fD' \subseteq \fD$ be any of the three classes of digraphs specified in Theorem \ref{theo_GschemeOnStrict}. Then, for all $R_1, R_2$ specified as in the theorem,
\begin{align*}
R_1 \strG S_1 \mytext{ and } R_2 \strG S_2 
\quad \Rightarrow & \quad R_1 + R_2 \strG S_1 + S_2.
\end{align*}
for all $S_1, S_2 \in \fD$.
\end{corollary}
\BP Let $G \in \fD'$ with connectivity components $G'_1, \ldots, G'_L$. Then, for all $H_1, H_2 \in \fD$,
\begin{align*}
\S(G, H_1 + H_2) & \; \simeq \; 
\prod_{\ell=1}^L \Big( \S(G'_\ell, H_1) + \S(G'_\ell, H_2) \Big). \end{align*}
Now the statement is a consequence of Theorem \ref{theo_GschemeOnStrict}.

\EP

\subsection{The rearrangement method} \label{subsec_rearrangement}

\begin{figure}
\begin{center}
\includegraphics[trim = 70 640 245 80, clip]{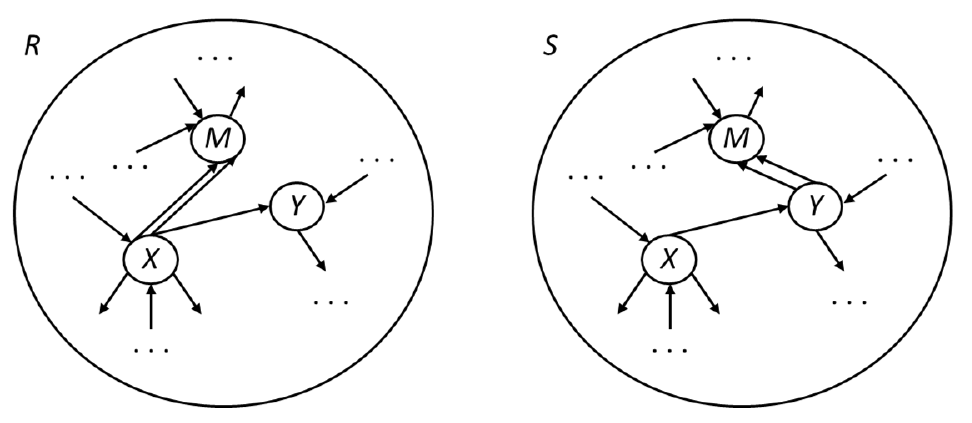}
\caption{\label{fig_abstrConstrMeth} The concept for the rearrangement of a digraph $R$. Explanations in text.}
\end{center}
\end{figure}

In this section, we present a method how to rearrange a digraph $R$ in such a way that $R \sqsubseteq_G S$ with respect to $\fD$ holds for the resulting digraph $S$. The method is applied on posets in Section \ref{subsec_posets}.

The concept of the rearrangement method is illustrated in Figure \ref{fig_abstrConstrMeth}. We have a digraph $R$ and sets $X, Y, M \subseteq V(R)$ with $X \cap M = \emptyset$, $Y \cap M = \emptyset$. We build a new digraph $S$ by replacing all arcs between $M$ and $X$ in $R$ by arcs between $M$ and $Y$. The approach is formalized as follows:

\begin{definition} \label{def_leq_s}
We agree on the following:
\begin{itemize}
\item $R = (Z, A(R))$ is a digraph.
\item We have disjoint subsets $X$ and $M$ of $Z$.
\item There is a subset $Y \subseteq Z$ with
\begin{align}
\label{bed_WB_empty}
M \cap Y & = \emptyset , \\
\label{bed_WBB_empty}
M \cap N_R(y) & = \emptyset \mytext{ for all } y \in Y,
\end{align}
and $ \beta : X \rightarrow Y $ is a mapping.
\item We define $S$ as the digraph with $V(S) = Z$ and 
\begin{align*}
A(S) & \equiv \; \; A_r \; \cup \; A_d \; \cup \;  A_u \\
\mytext{where} \; A_r & \equiv \; A(R) \setminus \left( ( M \times X ) \cup ( X \times M ) \right), \\
A_d & \equiv \;
\mysetdescr{ m \beta(x) }{ m x \in A(R) \cap ( M \times X )}, \\
A_u & \equiv \;
\mysetdescr{ \beta(x) m }{ x m \in A(R) \cap ( X \times M )}.
\end{align*}
\end{itemize}
\end{definition}

Every arc in $ A(R) \cap \left( ( M \mytimes X ) \cup ( X \mytimes M ) \right)$ is a proper arc due to $M \cap X = \emptyset$, and because of \eqref{bed_WB_empty}, $A_d$ and $A_u$ do not contain any loop. In consequence, the loops of $S$ are exactly the loops of $R$, and they all are contained in $A_r$. Furthermore, $A(R) \cap \left( A_d \cup A_u \right) = \emptyset$ according to \eqref{bed_WB_empty} and \eqref{bed_WBB_empty}.

\begin{definition} \label{def_Uxi}
For every $G \in \fD$ and every $\xi \in \S(G,R)$, we define the set $U_\xi$  by
\begin{equation*}
U_\xi \;  \equiv \; \mysetdescr{ v \in V(G) }{ \xiv \in X \mytext{and} \; \xi[ N_G(v) ] \cap M \not= \emptyset },
\end{equation*}
Correspondingly, for every $G \in \fD$ and every $\zeta \in \S(G,S)$, we define the set $U'_\zeta$ by
\begin{equation*}
U'_\zeta \;  \equiv \; \mysetdescr{ v \in V(G) }{ \zeta(v) \in Y \mytext{and} \; \zeta[ N_G(v) ] \cap M \not= \emptyset }.
\end{equation*}
\end{definition}

\begin{theorem} \label{theo_rhoxi}
For every $G \in \fD_r$ and every $\xi \in \S(G,R)$, we define the mapping $\rho_G(\xi) : V(G) \rightarrow Z$ by
\begin{eqnarray*}
\forall v \in V(G) \mytext{:} \; \rho_G(\xi)(v) & \equiv & 
\begin{cases}
\beta(\xiv), & \mytext{if } v \in U_\xi, \\
\xiv, & \mytext{otherwise}.
\end{cases}
\end{eqnarray*}
Assume additionally to the assumptions in Definition \ref{def_leq_s}, that $\beta : R \vert_X \rightarrow R \vert_Y $ is a bijective homomorphism and
\begin{align}
\label{bed_nbh}
\forall x \in X \mytext{:} \; \Nin_R(x) \setminus M \subseteq \Nin_R(  \beta(x) ) \; \mytext{and} \; \Nout_R(x) \setminus M \subseteq \Nout_R( \beta(x) ).
\end{align}
Then for every $G \in \fD_r$ and every $\xi \in \S(G,R)$,
\begin{eqnarray} \label{umkehrung_rhoxi}
\forall v \in V(G) \mytext{:} \; \xiv & = & 
\begin{cases}
\beta^{-1}(\rho_G(\xi)(v)), &  \mytext{if } v \in U'_{\rho_G(\xi)}, \\
\rho_G(\xi)(v), & \mytext{otherwise},
\end{cases}
\end{eqnarray}
and $\rho$ is a strong $\Gamma$-scheme from $R$ to $S$ with respect to $\fD$.
\end{theorem}
\BP Let $G \in \fD_r$ and $\xi \in \S(G,R)$ be selected. We write $U$ and $U'$ instead of $U_\xi$ and $U'_{\rho(\xi)}$, respectively. The reader should observe that for all $v \in Z$
\begin{align*}
\rhoxiv \not= \xiv & \; \Rightarrow \; v \in U, \; \xiv \in X,  \mytext{and} \rhoxiv = \beta( \xiv ), \\
\xiv \in M & \; \Rightarrow \; \rhoxiv = \xiv.
\end{align*}
To unburden the notation, we define 
\begin{equation*}
A_r^* \quad \equiv \quad A_r \setminus \Delta_R \quad = \quad A(R^*) \setminus (( M \times X) \cup (X \times M )).
\end{equation*}
Then $A(S^*) = A_r^* \cup A_d \cup A_u$.

{\em $\rho(\xi)$ is a strict homomorphism:} Let $vw \in A(G)$. If $vw$ is a loop in $G$, then $\rhoxiv \rhoxiv$ is a loop in $S$, because $\xi$ and $\beta \circ \xi \vert_U$ are both homomorphisms to $R$ and all loops in $R$ are loops in $S$, too. 

Assume that $vw$ is a proper arc in $G$. In the cases of $v, w \in U$ or $v, w \notin U$, we have $\rhoxiv \rhoxiw \in A^*_r$ due to the strictness of $\xi$ and $\beta \circ \xi \vert_U$, $X \cap M = \emptyset$ and \eqref{bed_WB_empty}. Let $v \in U$ and $w \notin U$. In the case of $\xiw \in M$, we have $\rhoxiv \rhoxiw \in A_u$, and in the case of $\xiw \notin M$, \eqref{bed_nbh} and \eqref{bed_WB_empty} yield $\rhoxiv \rhoxiw \in A^*_r$. Correspondingly, $v \notin U$, $w \in U$ yields $\rhoxiv \rhoxiw \in A^*_r \cup A_d$.

{\em Inversion formula \eqref{umkehrung_rhoxi}:} Looking at the definition of $\rho$, we realize that \eqref{umkehrung_rhoxi} holds if $U = U'$.

Let $v \in U$, thus $\xiv \in X$, $\rhoxiv \in Y$, and $\xi[ N_G(v) ] \cap M \not= \emptyset$. There exists a $w \in N_G(v)$ with $\xiw \in M$. We conclude $\xiw = \rhoxiw \in \rho(\xi)[ N_G(v) ]$, and $v \in U'$ is shown.

Let $v \in U'$, thus $\rhoxiv \in Y$ and $\rho(\xi)[ N_G(v) ] \cap M \not= \emptyset$. There exists a $w \in N_G(v)$ with $\rhoxiw \in M$. Then $\rhoxiw = \xiw$ according to \eqref{bed_WB_empty}, and due to \eqref{bed_WBB_empty}, the vertices $\rhoxiv$ and $\xiw$ are not adjacent in $R$. We conclude $\rhoxiv \not= \xiv$, hence $v \in U$.

{\em Final proof of $R \strG S$}: All together, we have proven that $\rho_G$ is a one-to-one mapping from $\S(G,R)$ to $\S(G,S)$ for all $G \in \fD_r$. Now apply Theorem \ref{theo_GschemeOnStrict}.

\EP

\section{Posets and undirected graphs} \label{sec_posets_graphs}

\subsection{Application on posets} \label{subsec_posets}

In this section, we show that for a poset $R$ fulfilling some additional conditions, the transitive hull $T$ of the digraph $S$ resulting from the rearrangement method is a poset with $R \strG T$ with respect to $\fD$.

Let $R \in \fD$ be as in Definition \ref{def_leq_s}, and assume that $R$ fulfills additionally the following conditions:
\begin{itemize}
\item $X$ is convex;
\item there is no walk in $R$ starting in $M$ and ending in $Y$, and no walk starting in $Y$ and ending in $M$.
\end{itemize}
The latter condition implies \eqref{bed_WBB_empty}, and it is not conflicting with \eqref{bed_nbh}. With $S$ defined as in Definition \ref{def_leq_s}, we get

\begin{lemma} \label{lemma_kettenregel}
Let $ z_0, \ldots , z_I $ be a walk in $S$. Define
\begin{align*}
K & \equiv \mysetdescr{ i \in \myNk{I} }{ z_{i-1} z_i \notin A_r }.
\end{align*}
Then $ \# K \leq 2$. If $K = \{ k, \ell \}$ with $k < \ell$, then there exist $m, n \in M$ and $x, y \in X$ with $mx \in A(R^*)$, $yn \in A(R^*)$ and
\begin{align*}
z_{k-1} z_k & = m \beta(x) \in A_d, \\
z_{\ell-1} z_\ell & = \beta(y) n \in A_u.
\end{align*}
\end{lemma}
\BP Let $\# K \geq 2$, and let $k < \ell$ be two {\em consecutive} indices in $K$ (i.e., $k, \ell \in K$ with $k < \ell$ and $ k < i < \ell \Rightarrow i \notin K$). There exist $x \in X$, $m \in M$ with $m x \in A(R^*)$, $z_{k-1} z_k = m \beta(x)$, or $x m \in A(R^*)$, $z_{k-1} z_k = \beta(x) m$, and there exist $y \in X$, $n \in M$ with $n y \in A(R^*)$, $ z_{\ell-1} z_\ell = n \beta(y)$, or $y n \in A(R^*)$, $z_{\ell-1} z_\ell = \beta(y) n$.

Let $ z_{k-1} z_k = m \beta(x)$, $z_{\ell-1} z_\ell = n \beta(y)$. With \eqref{bed_WB_empty}, we get $k < \ell - 1$, and $z_k, \ldots , z_{\ell-1}$ is a walk in $R$ starting with $\beta(x) \in Y$ and ending with $n \in M$ in contradiction to our assumption.

The case $z_{k-1} z_k = \beta(x) m$, $ z_{\ell-1} z_\ell = \beta(y) n$ is treated in the same way.

Now let $z_{k-1} z_k = \beta(x) m$ and $z_{\ell-1} z_\ell = n \beta(y)$. Then $x, z_k \ldots z_{\ell-1}, y$ is a walk in $R$, and the convexity of $X$ yields $m, n \in X$ in contradiction to $m, n \in M$.

The case $ z_{k-1} z_k = m \beta(x)$ and $z_{\ell-1} z_\ell = \beta(y) n$ remains as the only possible one. We learn from this result that for every pair $k' < \ell'$ of consecutive indices in $K$, we have $ z_{k'-1} z_{k'} \in M \times Y $ and $z_{\ell'-1} z_{\ell'} \in Y \times M $. Due to $ z_{k-1} z_k \in M \times Y$, $z_{\ell-1} z_\ell \in Y \times M$, and \eqref{bed_WB_empty}, the set $K$ cannot contain an index preceeding $k$ or an index following $\ell$, hence $K = \{ k, \ell \}$.

\EP

\begin{lemma} \label{lemma_leqS_is_partOrd}
If $R$ be an element of $\fTa$, then $S$ is an element of $\fTa$, too.
\end{lemma}
\BP Let $v \in Z$ and let $W = z_0, \ldots , z_I$ be a walk in $S$ with $z_0 = v = z_I$. As in Lemma \ref{lemma_kettenregel}, we define
\begin{align*}
K & \equiv \mysetdescr{ i \in \myNk{I} }{ ( z_{i-1}, z_i ) \notin \; A_r }.
\end{align*}
If $K = \emptyset$, then $W$ is a walk in $R$, hence $z_0 = \ldots = z_I$ due to $R \in \fTa$. In what follows, we show that indeed $K = \emptyset$. From Lemma \ref{lemma_kettenregel}, we know $\# K \leq 2$.

Let $K = \{ k \}$ with $k \in \myNk{I}$. There exist $m \in M, x \in X$ with $ z_{k-1} z_k = m \beta(x)$ or $ z_{k-1} z_k = \beta(x) m$. Assume $ z_{k-1} z_k = m \beta(x)$. The sequence $z_0, \ldots ,z_{k-1}$ is a walk in $R$ starting with $v$ and ending with $m$. In the case of $k < I$, the sequence $ z_k, \ldots , z_I $ is a walk in $R$ starting with $\beta(x)$ and ending with $v$, and in the case of $k = I$, we have $\beta(x) = v$. Therefore, in both cases, there exists a walk in $R$ starting with $\beta(x)$ and ending with $m$ in contradiction to our assumption. The proof for the case $ z_{k-1} z_k = \beta(x) m$ is analogous.

Now let $K = \{ k, \ell \}$ with $k < \ell$. According to Lemma \ref{lemma_kettenregel}, there exist $m, n \in M$ and $x, y \in X$ with $m x \in A(R^*)$, $ z_{k-1} z_k = m \beta(x)$ and $y n \in A(R^*)$, $ z_{\ell-1} z_\ell = \beta(y) n$. Now $ y, n, z_{\ell+1}, \ldots , z_I, z_1, \ldots, z_{k-2}, m, x$ is a walk in $R$, and the convexity of $X$ yields $m, n \in X$ in contradiction to $m, n \in M$.

\EP

\begin{corollary} \label{coro_rhoxi_posets}
Assume that $R$ is a poset and that the additional assumptions in Theorem \ref{theo_rhoxi} are fulfilled. Then the transitive hull $T$ of $S$ is a poset with $R \strG T$ with respect to $\fD$.
\end{corollary}
\BP Due to $\fP \subset \fTa$, we have $S \in \fTa$ according to Lemma \ref{lemma_leqS_is_partOrd}. Because all loops in $R$ are preserved, $S$ is reflexive, and $T$ is a poset. We have $R \strG S$ according to Theorem \ref{theo_rhoxi}, and the identity mapping of $Z$ is a one-to-one homomorphism from $S$ to $T$.

\EP

\begin{figure}
\begin{center}
\includegraphics[trim = 70 550 210 60, clip]{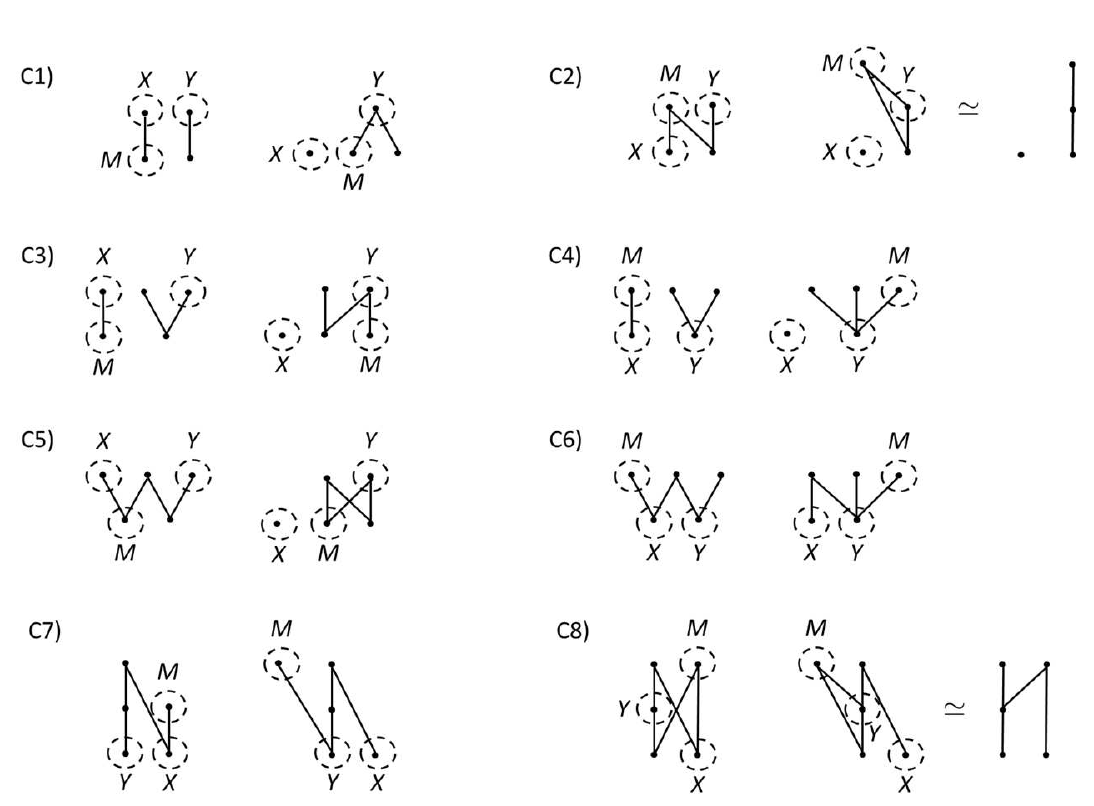}
\caption{\label{fig_TableConstr} Hasse-diagrams of eight examples for the application of the rearrangement method on posets. In all examples, $R$ is on the left and $T$ on the right. Explanations in text.}
\end{center}
\end{figure}

Figure \ref{fig_TableConstr} shows the Hasse-diagrams of eight pairs of posets $R$ and $T$. The respective one-element vertex sets $X$, $Y$, and $M$ in $R$ fulfill all conditions imposed in Definition \ref{def_leq_s}, Theorem \ref{theo_rhoxi}, and at the beginning of this section, and Corollary \ref{coro_rhoxi_posets} delivers $R \strG T$ in all eight cases.

For the first five pairs in the figure, we have $T = E + T'$ with a singleton $E$ and a connected poset $T'$. Theorem \ref{theo_GschemeOnStrict} sheds an interesting light on the roles of $E$ and $T'$ in $R \strG S$ in these five examples. For all connected digraphs $G$ with at least two points, $\# \S(G,R) \leq \# \S(G,E+T')$ is equivalent to $\# \S(G,R) \leq \# \S(G,T')$. It is thus $T'$ which bears the main burden of $R \strG T$, whereas the role of $E$ is to ensure $\# \S(E,R) \leq \# \S(E,T)$ by providing sufficiently many points in $T$.

\begin{figure}
\begin{center}
\includegraphics[trim = 60 640 370 70, clip]{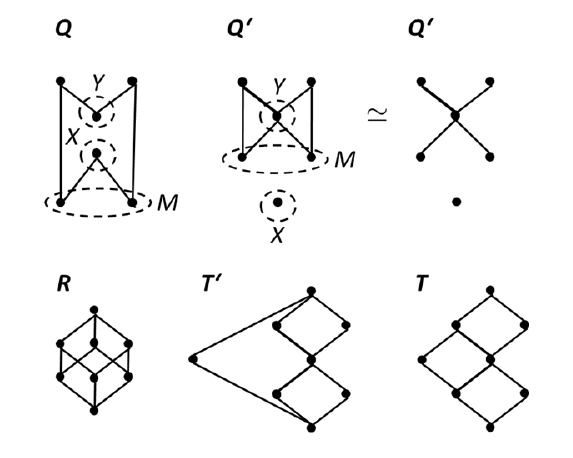}
\caption{\label{fig_Dedekind} Illustration for the proof of $R \strG T$ for the posets in Figure \ref{fig_Intro}(a). Explanations in text.}
\end{center}
\end{figure}

We have seen in Corollary \ref{calc_dirsum} that the relation $\strG$ is compatible with the direct sum of digraphs. It is also compatible with the ordinal sum of posets:

\begin{corollary} \label{calc_ordsum}
Let $R_1, R_2, S_1, S_2 \in \fP$. Then, with respect to $\fP$,
\begin{align*}
R_1 \sqsubseteq_\Gamma S_1 \mytext{ and } R_2 \sqsubseteq_\Gamma S_2 
\quad \Rightarrow & \quad R_1 \oplus R_2 \sqsubseteq_\Gamma S_1 \oplus S_2.
\end{align*}
\end{corollary}
\BP Let $P, Q_1, Q_2 \in \fP$, and let $\U(P)$ be the set of upsets of $P$. Then
\begin{align*}
\S(P, Q_1 \oplus Q_2) & \; \simeq \; \sum_{U \in \U(P)} \S(P \vert_{X \setminus U}, Q_1) \times \S(P \vert_U, Q_2).
\end{align*}
Now apply Theorem \ref{theo_GschemeOnStrict}.

\EP

Using the last two corollaries, we can provide a simple proof for $R \strG T$ with $R$ and $T$ as in Figure \ref{fig_Intro}(a). On top of Figure \ref{fig_Dedekind}, the Hasse-diagrams of two posets $Q$ and $Q'$ are shown for which Corollary \ref{coro_rhoxi_posets} yields $Q \strG Q'$. By applying Corollary \ref{calc_ordsum}, we get for the posets in the lower part of the figure 
\begin{equation*}
R \; \simeq A_1 \oplus Q \oplus A_1 \; \strG A_1 \oplus Q' \oplus A_1 \; \simeq T',
\end{equation*}
and we conclude $R \strG T$ because there exists a one-to-one homomorphism from $T'$ to $T$.

\subsection{Application on undirected graphs} \label{subsec_ugraphs}

It is commonly known that the category $\fS$ of {\em symmetric digraphs} and the category $\fU$ of {\em undirected graphs} are twins. For a digraph $G \in \fD$, we get the {\em underlying graph $u(G) \in \fU$} by replacing the arc set $A(G)$ by the {\em edge set} $E(u(G)) \equiv \mysetdescr{ \{ v, w \} }{ vw \in A(G) }$. In fact, $u : \fS \rightarrow \fU$ is a bijective covariant functor. Because homomorphisms are mappings between vertex sets, we have $u(\xi) = \xi$ for all $\xi \in \H(G,H)$, $G, H \in \fS$, and $\H_u( u(G), u(H) ) = \H( G, H )$, where $\H_u$ indicates the homomorphism set between {\em undirected} graphs in order to avoid confusion.

For $G \in \fU$, we denote by $G^*$ the undirected graph $G$ with loops (one-element-edges) removed, and we define the set of  {\em strict} homomorphisms between $G, H \in \fU$ as
\begin{align*}
\S_u(G,H) & \; \equiv \; \H_u(G,H) \cap \H_u(G^*,H^*),
\end{align*}
The operator $^*$ commutes with $u$, hence $\H(G^*,H^*) = \H_u(u(G)^*,u(H)^*)$ and $\S(G,H) = \S_u(u(S),u(H))$ for all $G, H \in \fS$.

We can transfer the definition of $\gxiv$ and $\Gxi$ without any modification into the world of undirected graphs, because the fundamental Definition \ref{def_connected} refers to the underlying graph. It is easily seen, that $\Gxi \in \fS$ holds for every homomorphism $\xi$ between symmetric digraphs and that also $\G$ commutes with $u$: $u( \Gxi ) = \myGxi{ u( \xi ) }$ for every homomorphism $\xi$ between symmetric digraphs. In particular, $\Gzet \in \fU$ for every homomorphism $\zeta$ between undirected graphs.

Therefore, we can define the relation $R \strG S$ for undirected graphs $R$ and $S$ by replacing digraphs by undirected graphs and $\H$ by $\H_u$ in Definition \ref{def_RS_scheme}. Moreover, we can transfer all our results to undirected graphs, as long as they are compatible with symmetry. That includes the complete Sections \ref{subsec_connectivity} and \ref{subsec_Def_HomSchemes}, and Section \ref{subsec_GschemeOnStrict} from Definition \ref{def_pixi_ioxi} until Lemma \ref{lemma_GDR}; the latter one reads for undirected graphs as follows:
\begin{lemma} \label{lemma_GDR_UG}
Let $R, S \in \fU$, let $\fU' \subseteq \fU$, and let
\begin{equation*}
\fG_{\fU'}(R) \quad \equiv \quad \mysetdescr{ \Gzet }{ \zeta \in \H_u(G,R), G \in \fU' }.
\end{equation*}
Then
\begin{align*}
\# \S_u(G,R) & \leq \# \S_u(G,S) \; \; \mytext{for all} \; G \in \fG_{\fU'}(R) \\
\mytext{implies} \quad \quad \quad \quad \quad \quad R &\strG S \; \; \mytext{with respect to} \; \fU'.
\end{align*}
\end{lemma}

Because the definition of the class $\fTa$ refers to antisymmetry, the results in connection with this class cannot be transferred to undirected graphs. But there is a class of undirected graphs for which similar results can be achieved. If we define
\begin{equation*}
\fC_o \; \equiv \; \mysetdescr{ G \in \fU }{ G^* \mytext{\em  does not contain a cycle of odd length} },
\end{equation*}
we have the following result corresponding to Lemma \ref{lemma_iota}:
\begin{lemma} \label{lemma_iota_UG}
For $G \in \fU$, $H \in \fC_o$ and $ \zeta \in \H_u(G,H)$, there exists no walk $\fc_0, \ldots , \fc_I$ in $\Gzet^*$ of odd length with $\iota_\zeta( \fc_0 ) = \iota_\zeta( \fc_I )$. In particular, $\Gzet \in \fC_o$.
\end{lemma}
(Start the proof with a walk $\fc_0, \ldots , \fc_I$ of odd length in $\Gzet^*$, proceed as in the original, and observe, that in the case of $\iota_\zeta( \fc_0 ) = \iota_\zeta( \fc_I )$, the walk $\zeta(v_0^+), \zeta(v_1^+), \ldots , \zeta(v_{I-1}^+), \zeta(v_I^-)$ is a cycle of odd length in $H^*$ because $\iota_\zeta$ is strict.)

Using this result and Lemma \ref{lemma_GDR_UG}, Theorem \ref{theo_GschemeOnStrict} becomes

\begin{theorem} \label{theo_GschemeOnStrict_UG}
Let $R \in \fU$ and $\fU' \subseteq \fU$. Then, for all $S \in \fU$, the equivalence
\begin{align*}
R & \strG S \; \; \mytext{with respect to } \fU'\\ \Leftrightarrow \quad \# \S_u(G,R) & \leq \# \S_u(G,S) \; \; \mytext{for all} \; G \in \fU',
\end{align*}
and the implication
\begin{align*}
\# \S_u(G,R) & \leq \# \S_u(G,S) \; \; \mytext{for all} \; G \in \fU' \\
\Rightarrow \quad \# \H_u(G,R) & \leq \# \H_u(G,S) \; \; \mytext{for all} \; G \in \fU'
\end{align*}
hold if
\begin{align*}
R & \in \fU' = \fU, \\
R & \in \fC_o \subseteq \fU' \subseteq \fU.
\end{align*}
\end{theorem}

Also the rearrangement method in Definition \ref{def_leq_s} can be transferred to undirected graphs. For $R \in \fU$ with the properties described in Definition \ref{def_leq_s}, define $S \in \fU$ by
\begin{align*}
V(S) & \equiv \; Z, \\
E(S) & \equiv \; \left( E(R) \setminus E_{M,X} \right) \; \cup \; E_b, \\
\mytext{where} \quad \quad E_{M,X} & \equiv \; \mysetdescr{ e \in E(R) }{ e \cap M \not= \emptyset \mytext{ and } e \cap X \not= \emptyset} \\
\mytext{and} \quad \quad \quad \; E_b & \equiv \;
\mysetdescr{ ( e \cap M ) \cup \beta[ e \cap X ] }{ e \in E_{M,X} }.
\end{align*}
Then the counterpart of Theorem \ref{theo_rhoxi} delivers $R \strG S$ if we replace \eqref{bed_nbh} by
\begin{align*}
\forall x \in X \mytext{:} & \; N_R(x) \setminus M \subseteq N_R( \beta(x) ). \end{align*}

\end{document}